\documentclass[11pt]{amsart}

\usepackage[square,compress,comma, numbers,sort]{natbib}
\usepackage[colorlinks=true, citecolor=red, linkcolor=blue]{hyperref}
\usepackage{amsfonts,mathtools}

\allowdisplaybreaks[4]

\usepackage{amsmath}
\usepackage{bbm}
\usepackage{amssymb}
\usepackage{mathtools}

\usepackage{amsmath, amssymb, amsfonts, amsthm}

\usepackage{color}
\definecolor{c20}{rgb}{0.,0.7,0.}
\definecolor{c30}{rgb}{0.,0.,1.}
\definecolor{c40}{rgb}{1,0.1,0.7}
\definecolor{c50}{rgb}{1,0,0}
\definecolor{c60}{rgb}{1,0.9,0.1}

\definecolor{c20}{rgb}{0.,0.7,0.}
\definecolor{c30}{rgb}{0.,0.,1.}
\definecolor{c40}{rgb}{1,0.1,0.7}
\definecolor{c50}{rgb}{1,0,0}
\definecolor{c60}{rgb}{1,0.9,0.1}

\newcommand{\ve}{\varepsilon}

\newcommand{\abs}[1]{\left\lvert #1 \right\rvert}

\newcommand{\E}[1]{\mathbb{E}\left\{ #1\right\}}

\newcommand{\pk}[1]{\mathbb{P} \left\{ #1 \right \} }
\newcommand{\PP}{\mathbb{P}}

\newcommand{\R}{\mathbb{R}}

\newcommand{\Z}{\mathbb{Z}}
\newcommand{\ups}{\Upsilon}

\newcommand{\inr}{\in \R}

\newcommand{\limit}[1]{\lim_{#1 \to   \infty}}

\newcommand{\BQN}{\begin{eqnarray}}
\newcommand{\EQN}{\end{eqnarray}}
\newcommand{\BQNY}{\begin{eqnarray*}}
\newcommand{\EQNY}{\end{eqnarray*}}

\newcommand{\BS}{\begin{sat}}
\newcommand{\ES}{\end{sat}}
\newcommand{\BT}{\begin{theo}}
\newcommand{\ET}{\end{theo}}
\newcommand{\BK}{\begin{korr}}
\newcommand{\EK}{\end{korr}}

\newcommand{\PH}{\overline{\Phi}}

\newcommand{\BD}{\begin{de}}
\newcommand{\ED}{\end{de}}

\newcommand{\BIT}{\begin{itemize}}
\newcommand{\EIT}{\end{itemize}}
\newcommand{\BDI}{\begin{description}}
\newcommand{\EDI}{\end{description}}

\newcommand{\BRM}{\begin{remarks}}
\newcommand{\ERM}{\end{remarks}}

\newcommand{\BEL}{\begin{lem}}
\newcommand{\EEL}{\end{lem}}

\def\bqny#1{\begin{eqnarray*} #1 \end{eqnarray*}}
\def\bqn#1{\begin{eqnarray} #1 \end{eqnarray}}

\newtheorem{theo}{Theorem}[section]
\newtheorem{sat}[theo]{Proposition}
\newtheorem{de}[theo]{Definition}
\newtheorem{lem}[theo]{Lemma}

\newtheorem{korr}[theo]{Corollary}
\newtheorem{remark}[theo]{Remark}
\newtheorem{remarks}[theo]{Remarks}

\newcommand{\netheo}[1]{{Theorem \ref{#1}}}

\newcommand{\prooftheo}[1]{ \textsc{\bf Proof of Theorem} \ref{#1}:}

\newcommand{\COM}[1]{}

\newcommand{\hp}{\widetilde{p}}
\newcommand{\hP}{\widetilde{P}}

\def\IF{\infty}

\newcommand{\QED}{\hfill $\Box$}

\def\ve{\varepsilon}

\def\IF{\infty}

\begin{document}

\title{Approximation of Ruin Probability and Ruin Time in Discrete Brownian Risk Models }

\author{{Grigori Jasnovidov}}
\address{Grigori Jasnovidov, Department of Actuarial Science, %\\Faculty of Business and Economics\\
	University of Lausanne,\\
	UNIL-Dorigny, 1015 Lausanne, Switzerland
}
\email{Grigori.Jasnovidov@unil.ch}

\bigskip

\date{\today}
 \maketitle

 {\bf Abstract:}  We analyze the classical Brownian risk models  discussing the approximation of ruin probabilities (classical, $\gamma$-reflected, Parisian and cumulative Parisian) for the case that ruin can occur only on specific discrete grids. A practical and natural grid of points is for instance $G(1)= \{0, 1, 2, \ldots\}$, which allows us to study the probability of the ruin on the first day, second day, and so one. For such a discrete setting, there are no explicit formulas for the ruin probabilities mentioned above.   In this contribution  we derive accurate approximations of ruin probabilities for uniform grids by letting the initial capital to grow to infinity.
 \\
 {\bf Key Words:}  Brownian motion; $\gamma$-reflected risk model; Parisian ruin probability; cumulative Parisian ruin;  ruin time approximation
\\
\\
 {\bf AMS Classification:} Primary 60G15; secondary 60G70

\section{Introduction}

The classical Brownian risk model of an insurance portfolio
\[R_u(t)= u+ ct -   B(t), \ \ \ \ \ t\ge 0,\]
with $B$ a standard Brownian motion,  the initial capital $u>0$ and
the premium rate $c>0$,
is a key benchmark model in risk theory; see e.g., \cite{ig:69}.
For any $u>0$ define the ruin time
$$ \tau(u)= \inf\{ t\ge 0: B(t)- c t> u  \}$$
and thus the  corresponding ruin probability is given by the well-known formula (see e.g., \cite{DeM15})
\bqn{ \label{contclasprob}
{ \psi}_\IF( u):=\pk{\tau(u)< \IF}= \pk{ \inf_{t\ge 0} R_u(t) < 0}= e^{- 2cu} .
}

In insurance practice however the ruin probability is relevant not on a continuous time scale, but on a discrete one, due to the operational time (which is discrete). %since ruin can happen at some given day, say.
For a given  discrete uniform grid $G(\delta)=\{0, \delta, 2 \delta , .... \} $ of step  $\delta>0$ we define the corresponding ruin probability by
 \BQN
 { \psi}_{\delta, \IF}(u):=\pk{ \inf_{t\in G(\delta)} R_u(t) < 0}= \pk{ \sup_{t\in G(\delta)} (  B(t)- ct)> u} .
 \EQN
For any $u>0$ it is not possible to calculate ${ \psi}_{\delta, \IF}(u)$ explicitly and no formulas are available for the distributional characteristics of the corresponding ruin time which we shall denote by $\tau_\delta(u)$.\\
A natural question when explicit formulas are lacking  is how can we approximate $ { \psi}_{\delta, \IF}(u)$  and $\tau_\delta(u)$ for large $u$? Also of interest is to know what is the role of $\delta$: does it influence the ruin probability in this classical risk model? The first question has been considered recently in \cite{Piterbargdiscrete} for fractional Brownian motion risk process.\\
When dealing with the Brownian risk model, both the independence of increments and the self-similarity property are crucial. In particular, those properties are the key to a rigorous and (relatively) simple proof.

Our first result presented next shows that the grid plays a role only with respect to the pre-factor specified by some constant. Specifically, that  constant is well-known in the extremes of Gaussian processes being the  Pickands
constant $\mathcal{H}_{2c^2 \delta }$, where
\bqn{\label{pcC} \mathcal{H}_\eta=	\frac{1}{\eta }\E{ \frac{\sup_{t\in \eta \mathbb{Z}} e^{ W(t)} }{  \sum_{ t\in \eta \mathbb{Z}} e^{ W(t) } }} =
		\frac{1}{\eta }\E{ \max_{ t\ge 0, t\in \eta \mathbb{Z}} e^{ W(t )}- \max_{ t\ge \eta, t\in \eta \mathbb{Z}} e^{ W(t )}}  \in (0,\IF)
		}
for any  $\eta>0$, with  $W(t)= \sqrt{2} B(t) - \abs{t}$. Note that   $\mathcal{H}_0=1$; the first formula in \eqref{pcC} is  derived in \cite{DiekerY}, whereas the second in \cite{KrzPMS}. In the following
 $\sim$ stands for asymptotic equivalence as $u\to \IF$.

\BT  \label{prop1} For any $\delta > 0$ we have
\bqn{  \label{appRT0}
\psi_{\delta,\IF} (u) \sim
\mathcal{H}_{2c^2\delta}  \psi_\IF(u), \quad u\to \IF
}
and further for any $s\inr$ with $\Phi$ being the distribution
function of a standard Gaussian random variable
\bqn{\label{appRT}
	\limit{u} \pk{ c^{3/2}  (\tau_\delta(u) - u /c)/\sqrt{u} \le s \Bigl \lvert \tau_\delta(u) < \IF} =\Phi(s).
}
\ET

We note that the above results hold for the continuous case too, where the grid $G(\delta)$ is substituted by $[0,\IF)$ i.e., $\delta=0$. For that case \eqref{appRT} follows from \cite{MR2462285}. The approximation in \eqref{appRT} shows that the ruin time is not affected by the density of the grid  (i.e., it is independent of $\delta$) and thus we conclude that the grid influences only the ruin probability. This is not the case for the ruin probability approximated in \eqref{appRT0}. For the Pickands constants we have, see e.g., \cite{DiekerY,MR3745392,KDEH1}
$$ \mathcal{H}_{2c^2\delta} \le  \mathcal{H}_{0}=1 = \lim_{\delta  \downarrow 0} \mathcal{H}_{2c^2\delta}  =
\lim_{c  \downarrow 0} \mathcal{H}_{2c^2\delta}.$$
In particular we see that via self-similarity in the Brownian risk
model the role of the grid is coupled with the premium rate $c>0$. \\
The objective of Section 2 is to explain in detail the main ideas and  techniques adequate for the classical Brownian risk model.
Section 3 discusses  the ruin probability for the $\gamma$-reflected
Brownian risk model, see also \cite{MR3091101,MR3256222,MR3379034,KEP20171}.
The approximation of Parisian ruin
(see \cite{LCPalowski,MR3414985,MR3457055})
and  cumulative Parisian ruin
(see \cite{KrzyszSjT,DebZbiXia,DHLanpengRolski})
is the topic of Section 4.  Our findings show that also for
those ruin probabilities, the influence of the grid, i.e., the
choice of $\delta$ concerns only the leading constant in the asymptotic
expansion being further coupled with the premium rate.
Given the technical nature of several proofs, we shall relegate  them to Section 5, which is followed by an Appendix containing auxiliary calculations.

\section{Approximation techniques for Brownian Risk Model}
Both the independence of increments and the self-similarity property of the Brownian motion render  the Brownian risk model very tractable.
In order to approximate $\psi_{\delta,\IF}(u)$ for given $\delta>0$  we start with the following lower  bound
$$ \psi_{\delta,\IF}(u) = \pk{ \exists t\in G(\delta) : B(t)> u+ ct} \ge \pk{ B(u t_u)>  u(1+ ct_u)} $$
valid for $t_u$ such that $ut_u\in G(\delta)$ for all $u$ large. It is clear that such $t_u$ exists  and moreover
\bqn{
t_u= \frac{1} c+ \frac {\theta_u}  u \in G(\delta)
}
holds for some  $\theta_u\in [0,\delta)$ and   all large $u$.
Consequently,
by the well-known Mill's ratio asymptotics
$\overline{\Phi}(u) \sim  \varphi(u)/u $ as $u\to \IF$
we obtain for all
large $u$  and some positive constant $C$
 \bqn{\label{low}
 	  \psi_{\delta,\IF}(u)  \ge \overline{\Phi}(   \sqrt{ u/t_u} (1+ c t_u) )
 \ge   \frac{C}{\sqrt{ u}} e^{- 2 cu},
}
where $\PH = 1-\Phi$ and $\varphi = \Phi'.$
Although the lower bound above is not precise enough, it is  useful to localize a short interval around
$$t_0:=1/c$$
 that will lead eventually to  the exact approximation  of the ruin probability.  Indeed, we have with
$$T_u^{\pm }= u (t_0 \pm u^{-1/2}\ln u), \ \ \ \ Z(t) = B(t)-ct$$
for all large  $u$ and any $C>0, \ p<0$ (the proof is given in
the Appendix)
\bqn{\label{up1}
	 \pk{ \sup_{t\notin [T_u^{-},T_u^{+}]} Z(t)> u} \le C u^p e^{-2 cu}.
}
Since for any $u>0$
\bqn{\label{ineq}
	 \pk{ \sup_{t \in   [T_u^{-},T_u^{+}] \cap \delta \mathbb{Z}} Z(t)> u}
	 \le   \psi_{\delta, \IF}(u) \le \pk{ \sup_{t \in   [T_u^{-},T_u^{+}] \cap \delta \mathbb{Z}} Z(t)> u}
 +  \pk{ \sup_{t\notin [T_u^{-},T_u^{+}]} Z(t)> u}
	}
by \eqref{low} and \eqref{up1} we obtain    that (set $\Delta_\delta(u)=[t_u- u^{-1/2} \ln u, t_u+ u^{-1/2} \ln u]\cap \frac{\delta}{u} \mathbb{Z}$)
$$  \psi_{\delta, \IF}(u)  \sim  \pk{ \sup_{t \in  [T_u^{-},T_u^{+}], t \in  \delta \mathbb{Z}} Z(t) > u}=
\pk{ \exists t \in \Delta_\delta(u)  : B(t)> \sqrt{u} (1+ ct)}=:P_\delta (u),    \quad u\to \IF ,
$$
where for the last equality we used the self-similarity property of Brownian motion.\\
In order to approximate $P_\delta(u)$ as $u\to \IF$
 a common approach is to partition $\Delta_\delta(u)$ in small intervals and use Bonferroni inequality in order to determine the main contribution to the asymptotics. This idea coupled with the continuous mapping theorem is essentially due to Piterbarg, see e.g.,  \cite{Pit96}. In this paper we use a modified approach in order to tackle some uniformity issues which arise in the approximations. In particular, we do not use continuous mapping theorem but rely instead on the independence of increments and self-similarity property of Brownian motion. We illustrate below briefly our approach. \\
We choose a partition $ \Delta_{j,S,u}, -N_u\le j\le N_u$ of $\Delta_\delta(u)$ depending on some constant $S>0$ as follows
%, where
% for any fixed $S,u$ positive
\bqn{\label{deltajsu}
\Delta_{j,S,u}= [t_u+j Su^{-1}, t_u+(j+1) Su^{-1}] \cap \frac{\delta}{u} \mathbb{Z}, \quad N_u=\lfloor S^{-1} \ln(u) \sqrt{u}\rfloor.
}
Here  $\lfloor\cdot\rfloor$ stands for  the ceiling function. The Bonferroni  inequality yields
\BQN\label{eq:thetaT}
p_1(S,u)\ge P_\delta(u) \ge p_1'(S,u)-p_2(S,u),
\EQN
where
\bqny{
p_1(S,u) & = & \sum_{j=-N_u-1}^{N_u}p_{j,S,u},\quad
p_1'(S,u)\ = \sum_{j=-N_u}^{N_u-1}p_{j,S,u}, \quad  p_2(S,u)  = \sum_{ -N_u-1\le j< i\le  N_u}p_{i,j;S,u},
}
with
\begin{equation*}
p_{j,S,u}=\pk{\exists_{t\in\Delta_{j,S,u} } B(t)>
	\sqrt{u}(1+ ct) }
\ \text{ and } \ p_{i,j;S,u}=\pk{\exists_{t\in \Delta_{i,S,u}}  B(t)
\! > \! \sqrt u (1+ ct) ,
 \exists_{t\in \Delta_{j,S,u}} B(t) \! > \! \sqrt u (1+ ct)    }.
\end{equation*}

As shown in \cite{Rolski17} [Eq.\ (43)] the term $p_2(S,u)$, also referred to as the double-sum term, is negligible compared with $p_1'(S,u)$ if we let $u\to \IF$ and then $S\to \IF$.\\
Moreover, $p_1(S,u)$ and $p_1'(S,u)$ are asymptotically equivalent with $P_\delta(u)$, i.e.,
$$  \limit{S}\limit{u} p_1(S,u)/p_1'(S,u) = \limit{S} \limit{u} p_1(S,u)/P_\delta(u)= 1.$$
The main question is therefore how to approximate $p_1(S,u)$?\\
In order to answer the above question we need to approximate each term $p_{j,S,u}$ as $u\to \IF$. Moreover, such approximation has to be  uniform for all $j$ satisfying $ -N_u \le j \le N_u$, which is a  subtle issue solved in this paper  by utilizing the independence of increments of Brownian motion and the self-similarity property;
see the proof of \netheo{prop1} in Section 5 and \cite{Rolski17} for similar ideas in the continuous time setting.

\section{$\gamma$-reflected Risk Model }
An interesting extension of the classical Brownian risk model is that of $\gamma$-reflected Brownian risk model introduced in \cite{MR2405335}.
	The  $\gamma$-reflected fractional Brownian motion risk model
	and its extensions are discussed in \cite{MR3091101, MR3256222,MR3379034, KEP20171}.
	In this section  we consider the approximation of the ruin probability over a discrete grid $G(\delta), \delta >0$ for the $\gamma$-reflected Brownian motion model.
Specifically,  for $\gamma \in (0,1)$ we define the risk model
	$$R_u^\gamma(t) =  u +ct- B(t)+ \gamma\inf\limits_{0 \le s \le t}(B(s)-cs), \quad t\ge 0.$$
For given $\delta>0$ we are interested in the ruin probability in discrete time, namely
	$$\Upsilon_{\gamma , \delta}(u) =
	\pk {\sup\limits_{t \in [0,\infty) \cap \delta \mathbb{Z}}
		\left(B(t)-ct-\gamma\inf\limits_{s \in [0,t]\cap \delta \Z}(B(s)-cs)\right)
		> u},$$
	which cannot be calculated explicitly.
The risk process $R_u^\gamma(t)$ is not Gaussian any more, however using
%the idea of \cite{MR3091101} we can transform the probability of interest. Namely,
 the independence of the increments of Brownian
motion and the self-similarity property,  for any $u>0$ we have with $t=k-l, \ s=l$
\bqny{
\Upsilon_{\gamma,\delta}(u) &=& \pk{\exists l\le k\in G(\delta): B(k)-k c-\gamma(B(l)-c l)>u}
\\ &=&
\pk{\exists l\le k\in G(\delta): (B(k)-B(l))+(1-\gamma)B(l)-c(k-\gamma l)>u}
\\&=& \pk{\exists l\le k\in G(\delta):
B(k-l)+(1-\gamma)B^*(l)-c(k-\gamma l)>u}
\\&=&
\pk{\exists t,s \in G(\delta): (B(t)-c t)+(1-\gamma)(B^*(s)-cs)>u}
\\&=&
\pk{\exists t,s \in G(\delta/u): \frac{B(t)+(1-\gamma)B^*(s)}{ct+(1-\gamma)cs+1}>\sqrt u},
%\\&:=&
%\pk{\exists t,s\in G(\delta/u): Z(t,s)>\sqrt u},
	}
	where $B^*$ is an independent copy of $B$. The above re-formulation shows that  the ruin probability concerns the supremum of the random field $Z$ given by
\bqn{ \label{Z}
	Z(t,s) =\frac{B(t)+(1-\gamma)B^*(s)}{ct+(1-\gamma)cs+1}, \quad s,t\ge 0.
}

From \cite{MR3091101} it follows, that for any $\eta,a>0$
$$\mathcal{P}^a_\eta := \E{\sup\limits_{t\in [0,\IF) \cap \eta\mathbb{Z}}
 e^{  \sqrt 2 B(t)-t(1+a)}}\in (0,\IF) .$$

Our next result gives the approximation of the above ruin probability as $u \to \IF$.
\begin{theo}\label{gammatheo} For any $\delta>0$ and any $\gamma\in (0,1)$
		\bqn{\label{claimgammatheo}
			\Upsilon_{\gamma,\delta}(u) \sim \mathcal{P}^{\frac{\gamma}{1-\gamma}}_{2c^2(1-\gamma)^2\delta}\mathcal{H}_{2c^2\delta}\psi_\IF(u), \ \ \ u \to \IF.
		}
\end{theo}

We note   that the basic properties of  discrete Piterbarg constants are discussed in \cite{longKrzys,KDEH1}.

\section{Parisian \& Cumulative Parisian Ruin}
\subsection{Parisian ruin}
In this section we expand our results to the Parisian ruin.  For the continuous time \cite{LCPalowski} gives an exact formula for the Parisian ruin probability.
Both finite and infinite Parisian ruin times for continuous setup of the
problem are dealt with  in \cite{MR3414985,MR3457055}.\\
Next, for given $\delta, T $ positive
(suppose for convenience that $T /\delta \in G(\delta)$) define the  Parisian ruin for the discrete grid $G(\delta)$ by
\bqny{
\mathcal{P}_\delta(u,T) = \pk{\inf\limits_{t\in [0,\infty) \cap \delta \mathbb{Z}}\sup\limits_{s\in [t,t+T] \cap \delta \mathbb{Z}}R_u(s)<0}
.}

Our next result shows again that the grid determines the asymptotic approximation via the constant $\mathcal{H}_{\eta, T}$ defined for $\eta, T$ positive by
 \bqn{\label{eqC}
 \mathcal{H}_{ \eta, T}=	\E{ \frac{\sup\limits_{t\in \eta  \mathbb{Z}}\inf\limits_{s\in[t, t+T] \cap   \eta \mathbb{Z}} e^{ \sqrt{2} B(s)- \abs{s}} }{ \eta \sum_{t\in \eta \mathbb{ Z} }  e^{ \sqrt{2} B(t)- \abs{t}}}} \in (0,\IF).
 }
 Note that if $T=0$, then $	\mathcal{H}_{ \eta, 0}$ equals the Pickands constant $	\mathcal{H}_{ \eta}$ defined in \eqref{pcC}. The corresponding constant for the continuous case is introduced in \cite{MR3414985}.

\begin{theo}  \label{theoparis}
For any $\delta,T>0$
 \bqn{  \label{asymparis}
 	\mathcal{P}_\delta(u,T) \sim
 \mathcal{H}_{2c^2\delta,2c^2T}  \psi_\IF(u), \ \ u \to \IF
 .}
\end{theo}

We see from the approximation above that the premium rate $c$ influences also the leading constant in the asymptotics.

\subsection{Cumulative Parisian ruin}
Cumulative  Parisian ruin for fractional Brownian motion risk model has been discussed recently in \cite{KrzyszSjT}.
As therein, adjusted for the discrete setup, we  define the cumulative Parisian ruin probability by
$$\mathcal{C}_\delta(u,k) = \pk{\# \{t \in G(\delta): B(t)-ct>u \}>k },$$
where $k$ is some non-negative integer  and the symbol $\#$
stands for  the number of the elements of a given set.
Note in passing that
$\mathcal{C}_\delta(u,0)=\psi_{\delta, \IF}(u)$.
Next, for $\eta>0$ define the constant
\bqny{
\mathcal{B}_{\eta}(k) = \limit{S}\frac{\mathcal{B}_{\eta}(S,k)}{S},
}
where for any $S>0$
\bqny{
\mathcal{B}_{\eta}(S,k) = \int\limits_\R
\pk{ \eta \sum_{ s \in [0,S] \cap \eta \mathbb{ Z}}
\mathbb{I}(\sqrt 2 B(s)-|s|+z>0)>k}
 e^{-z}dz,
}
with
 $\mathbb{I(\cdot)}$ denoting the indicator function. In view of
 \cite{DebZbiXia} $\mathcal{B}_{\eta}(k)$ is  positive and finite.
\begin{theo} \label{theocum}  For any non-negative integer $k$ we have as  $u \to \IF$
\bqn{\label{cumtheoprob}
\mathcal{C}_\delta(u,k)  \sim  \mathcal{B}_{2c^2\delta}(k)\psi_\IF(u)
.}
\end{theo}

\begin{remark}
i) Defining the ruin times corresponding to Parisian and cumulative Parisian ruin, it follows with similar arguments as in the proof of \netheo{prop1} that those can be approximated in the same way as \eqref{appRT}.\\
ii) If $k=0$, then the claim in \eqref{cumtheoprob} reduces to \eqref{appRT0}.
\end{remark}

\section{Proofs}
\prooftheo{prop1}
As mentioned in Section 2, the negligibility of the double-sum term follows by \cite{Rolski17}, hence  the claim in \eqref{appRT0}  follows thus by approximating $p_1(S,u)$ as $u\to \IF$. We show first the approximation of   $p_{j,S,u}$   as $u\to \IF$ uniformly for $ -N_u \le j \le N_u$.
Note that with $u=v^2$ and $\mathcal{N}$ being a standard Gaussian random variable we have the distributional representation based on the independence of increments of Brownian motion
$$B( c_{j,S,u}+ t/u) =\sqrt{ c_{j,S,v}} \mathcal{N}+ B(t)/v , \quad t\in [0,S], \ u>0, \ c_{j,S,v}= t_u+ j S v^{-2}.$$
Recall that $t_u \in G(\delta)$ is given by $t_u= 1/c+ \theta_u/u$ for some $\theta_u \in [0,\delta)$. It turns out that $\theta_u$ will not play any role in the final asymptotic approximation.  We have thus with $\varphi_{j,v}$
the probability density function of $ \sqrt{ c_{j,S,v}} \mathcal{N}$
\bqny{
	p_{j,S,u} &=&\pk{\exists_{t\in\Delta_{j,S,u} }: ( B(t) -  \sqrt{u} ct)>
		\sqrt{u}  }\\
	&=& \int_{\R} \pk{\exists_{t \in [0,S]\cap \delta \mathbb{Z}}  :( B(t)/v -  v c (c_{j,S,v}+ t/v^2)  > v-x
	\lvert \sqrt{ c_{j,S,v}} \mathcal{N} = x} \varphi_{j,v}(x) dx\\
	&=&\frac{1}{v} \int_{\R} \pk{\exists_{t \in [0,S]\cap \delta \mathbb{Z}}  :( B(t)/v -  v c (c_{j,S,v}+ t/v^2)  > v-(v-x/v)
} \varphi_{j,v}(v-x/v) dx\\
	&=&\frac{1}{v} \int_{\R} \pk{\exists_{t \in [0,S]\cap \delta \mathbb{Z}}  : Z(t)    >x+c c_{j,S,v}v^2
} \varphi_{j,v}(v-x/v) dx\\
	&=&\frac{1}{v} \int_{\R} \pk{\exists_{t \in [0,S]\cap \delta \mathbb{Z}}  : Z(t)   >x
} \varphi_{j,v}(v(1+cc_{j,S,v}) -  x/v) dx\\
	&=&\frac{ e^{ - v^2(1+cc_{j,S,v})^2/(2 c_{j,S,v})} }{v \sqrt{2 \pi c_{j,S,v} }} \int_{\R}
	w(x)\omega(j,S,x) dx,
}
where (recall $ Z(t)= B(t) -  ct,t\ge 0$)
\bqn{\label{omegajsx}
w(x)= \pk{\exists_{t \in [0,S]\cap \delta \mathbb{Z}}  : Z(t)  >x
} , \quad \omega(j,S,x)=e^{   x(1+cc_{j,S,v}) / c_{j,S,v  }-    x^2/(2 c_{j,S,v} v^2)}.
}
 Using Borell-TIS inequality (see e.g., \cite{LifBook}) we have
(proof is given in the Appendix)
\bqn{
\label{intlim}
\int_{\R} w(x) \omega(j,S,x) dx
 &=&
 \int_{-M}^M w(x) e^{   2cx} dx + A_{M,v},
}

where $A_{M,v}\to 0$ as $u\to \IF$ and then $M\to \IF$, uniformly for $- N_u \leq j \leq N_u$ and $S>0$.
By the monotone convergence theorem

\bqny{
\limit{M}\int_{-M}^M  w(x) e^{   2cx} dx
 =  \frac{1}{2c} \E{\sup_{t\in [0,S] \cap \delta
 \mathbb{Z}} e^{ 2 c B(t) - 2c^2 t}}
&=& \frac{1}{2c}\E{\sup\limits_{t \in [0,S] \cap  \delta \mathbb{Z}}e^{\sqrt 2 B(2c^2t)-2c^2t}}
\\& =&
\frac{1}{2c} \E{\sup_{t\in [0,2 c^2S] \cap 2c^2\delta \mathbb{Z}} e^{ \sqrt{2}  B(t) - t}}
.}
In a view of the definition of discrete Pickands constants, see e.g., \cite{DiekerY,SBK}
$$ \limit{S}\frac{1}{2c^2S}\E{\sup_{t\in [0,2 c^2S] \cap 2c^2\delta \mathbb{Z}} e^{ \sqrt{2}  B(t) - t}} =    \mathcal{H}_{2c^2\delta},
$$
with $\mathcal{H}_\eta$ defined in \eqref{pcC}. Consequently, the asymptotics of $p_1(S,u)$ as $u\to \IF$ and therefore also \eqref{appRT}  follow by calculating the limit as $u\to \IF, S\to \IF$ of
$$  K_{v,S}= e^{2 v^2 c}   c S \sum_{j=-N_u-1}^{N_u} \frac{ e^{ - v^2(1+cc_{j,S,v})^2/(2 c_{j,S,v})} }{v \sqrt{2 \pi c_{j,S,v} }}.$$

Setting
$$f(t)= (1+ c t)^2/2t = 1/(2t)+ c + c^2t/2, \quad f'(t)= (-1/t^2 +c^2)/2, \quad  f''(t)= 1/ t^3$$
we have that $f'(t_0)=0$ implying
$$f(t_0+ x)- f(t_0)= \frac{f''(t_0)}{2}x^2+O(x^3)$$
as $x \to 0$ with $f''(t_0)= c^3$. Consequently, as $u\to \IF$
\bqn{ K_{v,S} &\sim&
\frac c{ \sqrt{2 \pi t_0 }}       \frac{S} v
	\sum_{j=-N_u-1}^{N_u}
e^{ - \big( v^2 f(t_0+ (jS+ \theta_u)/v^2)  - v^2 f(t_0)\big) } \notag \\
	&\sim & \frac { c} { \sqrt{2 \pi t_0 }}   \frac{S} v
	\sum_{j=-N_u-1}^{N_u} e^{ - f''(t_0) ((jS+ \theta_u)^2/v^2) /2 } \notag  \\
	&\sim & \frac { c^{3/2}} { \sqrt{2 \pi  }}  \int_\R e^{ - f  ''(t_0) x^2/2}dx = 1 ,
	\label{ndihm}
}
where the last two steps follow with the same arguments as in the proof of (39) in \cite{Rolski17}. Finally, we have that as $u\to \IF$ and then $S \to \IF$
$$P_\delta(u) \sim K_{v,S}\mathcal{H}_{2c^2\delta}e^{-2uc}
 \sim \mathcal{H}_{2c^2\delta}e^{-2uc}.$$

We show next \eqref{appRT}. For any $u>0, s\in \R$ we have
\bqny{
	\pk{ \tau_{\delta}(u) -  u t_u  \le s \sqrt{u} \lvert \tau_{\delta}(u) < \IF } &=&
	\frac{1}{ \psi_{\delta, \IF}(u)} \pk{ \exists_{t \in [0, ut_u+  s \sqrt{u} ] \cap \delta \mathbb{Z} }:  Z(t)> u}.
}	
Considering the approximations of $p_{j,S,u}$ uniformly for all  $ -N_u \leq j \leq N_u'$ with $N_u'= \lfloor s  \sqrt{u} /S \rfloor $ we obtain as above
\bqny{
	\limit{u}	\frac{1}{ \psi_{\delta, \IF}(u)} \pk{ \exists_{t \in [0, ut_u+  s \sqrt{u} ] \cap \delta \mathbb{Z} }:  Z(t)> u}
&=&    \frac{ \int_{-\IF} ^s e^{ - f  ''(t_0) x^2/2}dx}{ \int_\R e^{ - f  ''(t_0) x^2/2}dx}
=  \frac{c^{3/2}}{\sqrt{2 \pi }} \int_{-\IF} ^s e^{ - c^3 x^2/2}dx
= \Phi(sc^{3/2}).
}
Hence
\bqny{
	\limit{u} \pk{ (\tau_\delta(u) - u t_u )/\sqrt{u} \le s \Bigl \lvert \tau_\delta(u) < \IF} =
	\Phi(sc^{3/2}) , \quad s \inr.
}

Since $\Phi$ is continuous, by Dini's theorem, the above convergence holds also substituting $s$ by $s_u$ such that  $\limit{u}s_u=s  \inr$.
Consequently, since  $\theta_u \in [0,\delta)$ we have also
\bqny{
	\limit{u} \pk{ c^{3/2}  (\tau_\delta(u) - u t_0  )/\sqrt{u} \le s \Bigl \lvert \tau_\delta(u) < \IF} =\Phi(s), \quad s \inr.
\ \ \ \ \ \ \ \ \ \  \ \ \ \ \ \ \ \ \ \ \ \  \ \ \ \hfill \Box}

\textbf{Proof of Theorem \ref{gammatheo}.}
Recall that $t_u = t_0+\theta_u/u = 1/c+\theta_u/u$ and denote $\beta=1-\gamma$.
We analyze the variance function $\sigma_Z^2$ of the process $Z(t,s)$. For any non-negative
$s,t$ we have
$$\sigma^2_Z(t,s) =  \frac{t+\beta^2s}{(ct+\beta cs+1)^2} =
\frac{t+\beta s}{(ct+\beta cs+1)^2}-\frac{\beta(1-\beta)s}{(ct+\beta cs+1)^2}=: A(t,s)-A^*(t,s).$$
Note that $A (t,s)$ depends only on $t+\beta s$ and achieves its global maxima on the line $t+\beta s = t_0=1/c$, while
$A^*(t,s)$ is negative for all $s>0$ and equals zero for $s=0$.
Hence
$(t,s) = (1/c,0) $ is the unique global maxima of  $\sigma_Z^2(t,s)$  and
$\sigma^2_Z(1/c,0) = \frac{1}{4c}.$ We define next %Denote below the important interval
$$\mathbb{D}_\delta(u) = \Bigl \{s,t \in G(\delta/u): (t,s)\in (-\frac{\ln u}{\sqrt u}+t_u,\frac{\ln u}{\sqrt u}+t_u)\times(0,\frac{\ln u }{\sqrt u}) \Bigr\}.$$
We have (proof see in the Appendix)
\bqn{\label{gammaasym}
	\ups_{\gamma,\delta}(u) \sim \pk{\exists t,s \in \mathbb{D}_\delta(u): Z(t,s)>\sqrt u}=:\zeta(u), \ \ u \to \IF.
}
%We take a similar partition as in the Brownian motion case, namely its expansion to the rectangle
% $[-\frac{\ln u}{\sqrt u}+t_u,\frac{\ln u}{\sqrt u}+t_u]_{\delta/u}
%\times [0,\frac{\ln u}{\sqrt u}]_{\delta/u}$.
Let $\Delta_{i,S,u}$ be as in \eqref{deltajsu} and set \bqny{
	p(i,j) = \pk{\exists t,s \in \Delta_{i,S,u} \times \Delta_{j,S,u}^*: (B(t)-\sqrt u c t)+\beta(B^*(s)-\sqrt u cs)>\sqrt u },\\
	p(i,j;i',j') = \PP\{\exists t,s \in \Delta_{i,S,u} \times \Delta_{j,S,u}^*: (B(t)-\sqrt u c t)+\beta(B^*(s)-\sqrt u cs)>\sqrt u, \\
	\exists t,s \in \Delta_{i',S,u} \times \Delta_{j',S,u}^*: (B(t)-\sqrt u c t)+\beta(B^*(s)-\sqrt u cs)>\sqrt u\}
}
 for $-N_u\le i \le N_u$, $0\leq j \le N_u$, fixed $S>0$ and
$$ \Delta_{j,S,u}^* = [\frac{jS}{u}, \frac{j(S+1)}{u}].$$

By Bonferroni inequality
$$
\sum\limits_{0\le j\le N_u-1,-N_u\le i \le N_u-1}p(i,j)-
\sum\limits_{0\le j,j'\le N_u,-N_u-1\le i,i' \le N_u, (i,i')\neq(j,j')}p(i,j;i',j')
\le \zeta(u)
\leq \sum\limits_{0\le j\le N_u,-N_u-1\le i \le N_u}p(i,j).
$$
The term
$$\sum\limits_{0\le j,j'\le N_u,-N_u-1\le i,i' \le N_u, (i,i')\neq(j,j')}p(i,j;i',j')$$
is negligible by the proof of Theorem 2.1, Eq. $[14]$ in \cite{MR3091101} and consequently

$$\zeta(u) \sim \sum\limits_{0\le j\le N_u,-N_u\le i \le N_u}p(i,j), \ \ u\to \IF.$$
Next we approximate $p(i,j)$ uniformly.
Recall, that $v^2 =u$, $c_{i,S,v}=t_u+\frac{iS}{v^2}$,
$\varphi_{i,v}$ is the density function of $\sqrt {c_{i,S,v}} \mathcal{N}$
and set
$G_j = [jS,(j+1)S]\cap \delta\Z$. Denote also $[a,b]_\alpha = [a,b]\cap\alpha\Z$
for any real numbers $a\le b$ and $\alpha>0$.
We have
\bqn{\label{p_ijexactformula}
p(i,j) &=&
\PP\{\exists (t,s)\in \Delta_{i,S,u}\times\Delta_{j,S,u}^*:
B(t)-B(c_{i,S,v})-c \sqrt u (t-c_{i,S,v})+B(c_{i,S,v})-\sqrt u cc_{i,S,v}
\notag\\ & \ &
+ \beta(B^*(s)-\sqrt u cs)>\sqrt u\}
\notag\\ &=&
\int\limits_{\R}\pk{\exists (t,s)\in [0,S/u]_{\delta/u}\times \Delta_{j,S,u}^*:B(t)-\sqrt u ct-\sqrt u cc_{i,S,v}
+\beta(B^*(s)-\sqrt u cs)>\sqrt u -x
}
\varphi_{i,v}(x)dx
\notag\\ &=&
\int\limits_{\R}
\pk{\exists (t,s) \in [0,S]_\delta\times G_j : \frac{B(t)}{v}-vc(c_{i,S,v}+\frac{t}{v^2})
	+\beta(\frac{B^*(s)}{v}-\frac{cs}{v})>v-x}
\varphi_{i,v}(x)dx
\notag\\ &=&
\frac{1}{v}\int\limits_{\R}
\pk{\exists (t,s) \in [0,S]_\delta\times G_j : \frac{B(t)}{v}-vc(c_{i,S,v}+\frac{t}{v^2})
	+\beta(\frac{B^*(s)}{v}-\frac{cs}{v})>v-(v-\frac{x}{v})}
\varphi_{i,v}(v-\frac{x}{v})dx
\notag\\ &=&
\frac{1}{v}\int\limits_{\R}
\pk{\exists (t,s) \in [0,S]_\delta\times G_j : B(t)-ct
	+\beta(B^*(s)-cs)>x+v^2c c_{i,S,v}}
\varphi_{i,v}(v-\frac{x}{v})dx
\notag\\  &=&
	\frac{1}{v}\int\limits_{\R}
\pk{\exists (t,s) \in [0,S]_\delta\times G_j : B(t)-ct
	+\beta(B^*(s)-cs)>x}
\varphi_{i,v}(v(1+cc_{i,S,v})-\frac{x}{v})dx
\notag\\  &=&
	\frac{e^{-\frac{v^2(1+cc_{i,S,v})^2}{2c_{i,S,v}}}}{v\sqrt{2\pi c_{i,S,v}}}
	\int\limits_{\R}
	W_j(x)\omega(i,S,x)dx,\notag
}
where
\bqny{
	W_j(x) = \pk{\exists (t,s) \in [0,S]_\delta\times G_j : B(t)-ct
		+\beta(B^*(s)-cs)>x}
}
and $\omega(i,S,x)$ is defined in \eqref{omegajsx}.
By Borell-TIS inequality for all $|i|,|j|\leq N_u$
(proof is in the Appendix)
\bqn{\label{intgammasim}
	\int\limits_{\R}
	W_j(x)\omega(i,S,x)dx
	\sim
	\int\limits_{\R}
	W_j(x)e^{2cx}dx, \ \ u \to \IF.
}

Next we have with $G_j^* = [2jc^2S,2(j+1)c^2S]\cap 2c^2\delta \Z$
\bqn{\label{finint}
	\int\limits_{\R}
	W_j(x)e^{2cx}dx &=&
	\frac{1}{2c}\int\limits_\R
	\pk{\sup\limits_{(t,s) \in [0,S]_\delta\times G_j} \left(2cB(t)-2c^2t+\beta(2cB^*(s)-2c^2s)\right)>2cx}e^{2cx}d(2cx)
\notag\\ &=&
\frac{1}{2c}\int\limits_\R
\pk{\sup\limits_{(t,s) \in [0,S]_\delta\times G_j} \left(\sqrt 2
B(2c^2t)-2c^2t+\beta(\sqrt 2B^*(2c^2s)-2c^2s)\right)>x}e^xdx
\notag\\&=&
\frac{1}{2c}\E{\sup\limits_{(t,s)\in [0,2c^2S]_{2c^2\delta}\times G_j^*}
\exp\left(\sqrt 2 B(t)-t+\beta(\sqrt 2 B^*(s)-s)\right)}
\notag\\&=&
\frac{1}{2c}	\E{\sup\limits_{t\in [0,2c^2S]_{2c^2\delta}}
e^{\sqrt 2 B(t)-t}}
\E{\sup\limits_{s \in G_j^*}e^{ \beta(\sqrt 2 B^*(s)-s)}}
.}
By %\eqref{p_ijexactformula},
\eqref{intgammasim} combined with the line above we write
\bqn{\label{asymgamma}
	\zeta(u) &\sim&
	\frac{1}{2c}
	\E{\sup\limits_{t\in [0,2c^2S]_{2c^2\delta}}
		e^{\sqrt 2 B(t)-t}} \notag
	\\ &\times&
	\sum\limits_{0\le j\le N_u}
	\E{\sup\limits_{s \in G_j^*}
		e^{\beta(\sqrt 2 B^*(s)-s)}}
	\sum\limits_{-N_u\le i \le N_u}
	\frac{e^{-\frac{v^2(1+cc_{i,S,v})^2}{2c_{i,S,v}}}}{v\sqrt{2\pi c_{i,S,v}}}, \ \ \ u \to \IF
	.}
As was shown in the proof of Theorem \ref{prop1} as $u \to \IF$ and then $S \to \IF$
\bqn{\label{mainsum}
	\frac{1}{2c}
	\E{\sup\limits_{t\in [0,2c^2S]_{2c^2\delta}}
		e^{\sqrt 2 B(t)-t}}
	\sum\limits_{-N_u\le i \le N_u}
	\frac{e^{-\frac{v^2(1+cc_{i,S,v})^2}{2c_{i,S,v}}}}{v\sqrt{2\pi c_{i,S,v}}}\sim \mathcal{H}_{2c^2\delta}e^{-2cu}.
}
We have as $S \to \IF$
(proof of the first line below is in the Appendix)
\bqn{\label{ratiosum}
		\sum\limits_{0\le j\le N_u}
	\E{\sup\limits_{s \in G_j^*}e^{\beta(\sqrt 2 B^*(s)-s)}}
	&\sim&
	\E{\sup\limits_{s \in [0,2c^2S]_{2c^2\delta}}
		e^{\beta(\sqrt 2 B(s)-s)}}
	\\&=&
	\int\limits_\R\pk{\sup\limits_{s\in [0,2c^2S]_{2c^2\delta}}\left(\sqrt 2 B(s\beta^2)-\frac{s\beta^2}{\beta}>x\right)}e^xdx
	\notag\\ &=&
	\int\limits_\R\pk{
\sup\limits_{s\in [0,2c^2\beta^2S]_{2c^2\beta^2\delta}}\left(
\sqrt 2 B(s)-
		s(1+\frac{1-\beta}{\beta})\right)>x}e^xdx
	\notag\\ &\rightarrow&
	\mathcal{P}^{\frac{1-\beta}{\beta}}_{2c^2\beta^2\delta}
	=
	\mathcal{P}^{\frac{\gamma}{1-\gamma}}_{2c^2(1-\gamma)^2\delta}
\in (0,\IF)
.\notag
}
Combining the statement above with \eqref{asymgamma} and \eqref{mainsum} we conclude
\bqny{
	\zeta(u) \sim
	\mathcal{P}^{\frac{\gamma}{1-\gamma}}_{2c^2(1-\gamma)^2\delta}
	\mathcal{H}_{2c^2\delta}e^{-2cu}, \ \ u \to \IF
}
and hence by \eqref{gammaasym} the claim follows.
\QED
\\
\\
\textbf{Proof of Theorem \ref{theoparis}.} The proof is similar to that of \netheo{prop1} and we use similar notation as therein.
We have by \eqref{up1}
\bqn{\label{hatp}
\mathcal{P}_\delta(u,T)
\sim
\pk{\sup\limits_{t \in [T_u^{-},T_u^+] \cap \delta \mathbb{Z}} \inf\limits_{s\in[t,t+T] \cap \delta \mathbb{Z}} Z(s)>u} =: \hP_\delta(u), \quad u\to \IF
}
if we show that $\hP_\delta(u) \ge Ce^{-2cu}$.
By the self-similarity of Brownian motion
\bqny{
\hP_\delta(u)  &=&
\pk{
 \exists t \in [T_u^{-},T_u^+] \cap \delta \mathbb{Z} : \forall
 s\in[t,t+T] \cap \delta \mathbb{Z}  Z(s)>u)} \\
 &=&
\pk{
 \exists t \in [\frac{T_u^{-}}{u},\frac{T_u^+}{u}] \cap \delta \mathbb{Z}:
  \forall  s\in[t,t+\frac{T}{u}] \cap \frac{\delta}{u} \mathbb{Z} B(s)>\sqrt u(1+cs)}
.}
We choose the same partition ${\Delta_{j,S,u}}, \ -N_u\le j\le N_u$ of the interval $\Delta_{u} = [\frac{T_u^-}{u},\frac{T_u^+}{u}]$ as in the proof of Theorem \ref{prop1}.
The Bonferroni inequality yields
\BQN
\hp_1(S,u)\ge \hP_\delta(u) \ge \hp_1'(S,u)-\hp_2(S,u),
\EQN
where
\begin{equation*}
\hp_1(S,u)  =  \sum_{j=-N_u-1}^{N_u}\hp_{j,S,u},\ \ \ \
\hp_1'(S,u)  =  \sum_{j=-N_u}^{N_u-1}\hp_{j,S,u},\ \ \ \
\hp_2(S,u)  =  \sum_{ -N_u-1\le j< i\le  N_u}\hp_{i,j;S,u},
\end{equation*}
with
\begin{equation*}
\hp_{j,S,u} = \pk{\sup\limits_{t\in \Delta_{j,S,u}}\inf\limits_{s\in [t,t+\frac{T}{u}] \cap \frac{\delta}{u} \mathbb{Z}}
	Z(s)>u}
\end{equation*}
and
$$
\hp_{i,j;S,u}=\pk{ \sup\limits_{t\in \Delta_{i,S,u}}\inf\limits_{s\in [t,t+\frac{T}{u}]\cap\frac{\delta}{u} \mathbb{Z}}Z(s)>u,
\sup\limits_{t\in \Delta_{j,S,u}}\inf\limits_{s\in [t,t+\frac{T}{u}] \cap \frac{\delta}{u} \mathbb{Z}} Z(s)>u   }.
$$
Clearly,
$\hp_{i,j;S,u} \leq p_{i,j;S,u}$ and hence
$$\hp_2(S,u) \ \leq \ p_2(S,u).$$
Thus, if we show that $\hp_1(S,u) \sim C_1e^{-2cu}$ we conclude
that $\hp_2(S, u)$ is negligible. We approximate each summand
in $\hp_1(S,u)$ uniformly. As in the proof of \netheo{prop1} we obtain
\bqny{
\hp_{j,S,u}
&=&
\frac{ e^{ - v^2(1+cc_{j,S,v})^2/(2 c_{j,S,v})} }{v \sqrt{2 \pi c_{j,S,v} }}
 \int_{\R} w(T,x) \omega(j,S,x)
dx
,}
where
$$w(T,x)=\pk{\sup\limits_{t\in [0,S] \cap \delta \mathbb{Z}}\inf\limits_{s\in [t,t+T] \cap \delta \mathbb{Z}} Z(s)>x
} $$
and $\omega(j,S,x)$ is defined in \eqref{omegajsx}.
By  Borell-TIS inequality (similarly the proof of \eqref{intlim})
it follows that
\bqny{%\label{intlimparis}
 \int_{\R} w(T,x) \omega(j,S,x)dx \to
 \int_\R w(T,x) e^{   2cx} dx, \ \ \ \ u \to \IF.
}
Next we have
\bqny{
\int_\R w(T,x)e^{   2cx} dx =  \frac{1}{2c} \E{\sup\limits_{t\in [0,S] \cap \delta \mathbb{Z}}
\inf\limits_{s\in [t,t+T] \cap \delta \mathbb{Z}}
e^{ 2 c B(s) - 2c^2 s}}
 &=&  \frac{1}{2c}\E{\sup\limits_{t \in [0,S] \cap  \delta \mathbb{Z}}
\inf\limits_{s\in [t,t+T] \cap \delta \mathbb{Z}}
e^{\sqrt 2 B(2c^2s)-2c^2s}}
\\
&=&
\frac{1}{2c} \E{\sup\limits_{t\in [0,2 c^2S] \cap 2c^2\delta \mathbb{Z}}
\inf\limits_{s\in [t,t+2c^2T] \cap 2c^2\delta \mathbb{Z}}
e^{ \sqrt{2}  B(s) - s}}
.}

It follows with similar arguments as in \cite{SBK} that as $S\to \IF$
\bqn{\label{pickparis}
	\limit{S} \frac{1}{2c^2 S}\E{ \sup\limits_{t \in [0, 2c^2S] \cap  2 c^2 \delta\mathbb{Z}}\inf\limits_{s\in [t,t+2c^2T] \cap  2c^2\delta \mathbb{Z}} e^{\sqrt2 B(t)-t}} = \mathcal{H}_{2 c^2\delta, 2Tc^2}\in (0,\IF),
}
where the constant $\mathcal{H}_{2c^2\delta ,2Tc^2}$
is given by \eqref{eqC}. Hence by \eqref{ndihm} we have
$$\hP_\delta(u) \sim \mathcal{H}_{2c^2\delta, 2Tc^2}e^{-2cu}, \ \ \ u \to \IF$$
and \eqref{hatp} holds, establishing the claim.
\QED
\\

\textbf{Proof of Theorem \ref{theocum}.}
We use below the same notation as in the previous proofs.
By \eqref{up1} we have
\bqn{\label{peterb}
\mathcal{C}_\delta(u,k) \sim
\pk{\# \{t \in G(\delta) \cap  (T_u^-,T_u^+): Z(t)>u \}>k} =:  \hat{\psi}_k^\delta(u), \ \ \ u \to \IF
}
if we show that $\hat{\psi}_k^\delta(u) \ge Ce^{-2cu}$.
Using the self-similarity of Brownian motion for any $u>0$
$$\hat{\psi}_k^\delta(u) = \pk{\# \{t \in (-\frac{\ln u}{\sqrt u}+t_0,t_0+\frac{\ln u}{\sqrt u}) \cap  G(\frac{\delta}{u}): \frac{B(t)}{ct+1}> \sqrt{u} \}>k}. $$

Letting  $$A_{j,u} := \# \{t \in \Delta_{j,S,u}: Z(t)>u \}$$
we have using the idea from \cite{DHLanpengRolski}
\bqny{
\hat{\psi}_\delta^k(u)  &\leq&
\pk{\sum\limits_{j = -N_u-1}^{N_u} A_{j,u}>k}
\\&=&
\pk{\sum\limits_{j = -N_u-1}^{N_u} A_{j,u}>k,
\text{\{there exists only one $j$ such that }A_{j,u}>0 \}}
\\ & \ & + \ \
\pk{\sum\limits_{j = -N_u-1}^{N_u} A_{j,u}>k,\text{\{there exists $i \neq j$ such that $A_{i,u}>0$ and $A_{j,u}>0$ \}}}
\\ &=:&  p_{1,k}(u)+\Pi_0(u)
.}
On the other hand
\bqny{
\hat{\psi}_k^\delta(u) &\geq&
\pk{\sum\limits_{j = -N_u}^{N_u-1} A_{j,u}>k}\\
& = &
\pk{\sum\limits_{j = -N_u}^{N_u-1} A_{j,u}>k,\text{\{there exists only one
		{$j$} such that   $A_{j,u}>0$ \}}}
\\ & \ & + \ \
\pk{\sum\limits_{j = -N_u}^{N_u-1} A_{j,u}>k,\text{\{there exists $i \neq j$ such that  $A_{i,u}>0$ and $A_{j,u}>0$ \}}}
\\
&=:&
p_{2,k}(u)+\Pi'_0(u).
}

Notice, that $\Pi_0(u)$ and $\Pi_0'(u)$ are less than the double-sum term in \netheo{prop1}. They are negligible if we prove that $p_{2,k}(u)\sim p_{1,k}(u) \ge Ce^{-2cu}$ as $u \to \IF$ for some $C>0$.
We have
\bqny{
p_{1,k}(u)
 &=&
\sum\limits_{j = -N_u-1}^{N_u} \Big(
\pk{A_{j,u}>k}-\pk{A_{j,u}>k, \exists i\neq j: A_{i,u}>0}\Big)
\\ &=& \sum\limits_{j = -N_u-1}^{N_u} \pk{A_{j,u}>k} - \sum\limits_{j = -N_u-1}^{N_u} \pk{A_{j,u}>k, \exists i\neq j: A_{i,u}>0}.
}
The last summand is less than the double-sum term in \netheo{prop1}
and is negligible.
Thus we need to compute the asymptotic of
\bqn{\label{sumcum}
Q_{\delta,k}(u) := \sum\limits_{j = -N_u}^{N_u} \pk{A_{j,u}>k}
.}

With similar arguments as in the proof of \netheo{prop1}
\bqny{
\pk{A_{j,u}>k}  = \pk{ \# \{t \in \Delta_{j,S,u}: Z(t)>u \} > k}
 =
\frac{ e^{ - v^2(1+cc_{j,S,v})^2/(2 c_{j,S,v})} }
{v \sqrt{2 \pi c_{j,S,v} }} \int_{\R}w_k(x)
\omega(j,S,x) dx,
}
where $\omega(j,S,x)$ is defined in \eqref{omegajsx} and
 $$w_k(x)=
\pk{\#\{t \in [0,S]\cap \delta \mathbb{Z}: Z(t) >x\}>k}.$$

Similarly to the proof of \eqref{intlim} we have
\bqny{\label{cumlim}
\int_{\R} w_k(x)\omega(j,S,x)dx
 \to
\int_\R w_k(x)e^{ 2cx} dx, \ \ \ \ u \to \IF.
}
Next
\bqny{
\int\limits_\R w_k(x) e^{ 2cx}dx
=
\frac{1}{2c}\int\limits_\R \pk{\#\{t \in [0,2c^2S]\cap 2c^2 \delta \mathbb{Z}: \sqrt 2 B(t)-t >x\}>k}e^xdx
.}
As shown in \cite{DebZbiXia}
\bqny{ \limit{S}
\frac{1}{2c^2S}\int\limits_\R \pk{\#\{t \in [0,2c^2S]\cap 2c^2 \delta \mathbb{Z}: \sqrt 2 B(t)-t >x\}>k}e^xdx
= \mathcal{B}_{2c^2\delta}(k) \in (0,\IF)
.}

Consequently, by \eqref{ndihm}
as $u \to \IF$ and then $S \to \IF$
\bqny{
Q_{\delta,k}(u) \sim  cS\mathcal{B}_{2c^2\delta}(k)
\sum\limits_{j=-N_u}^{N_u}
\frac{ e^{ - v^2(1+cc_{j,S,v})^2/(2 c_{j,S,v})} }{v \sqrt{2 \pi c_{j,S,v} }}
& =&  \mathcal{B}_{2c^2\delta}(k)e^{-2cu}K_{v,S}  \sim
\mathcal{B}_{2c^2\delta}(k)e^{-2cu}.
}
Since $\mathcal{B}_{2c^2\delta}(k)\in (0,\IF)$, then  $Q_{\delta,k}(u) \sim \hat{\psi}_k^\delta(u)$ as $u \to \IF$ implying
\bqny{
\hat{\psi}_k^\delta(u) \sim e^{-2cu}\mathcal{B}_{2c^2\delta}(k),
 \quad u\to \IF.
\ \ \ \ \ \ \ \ \ \ \ \ \ \ \ \ \ \ \ \
 \ \ \ \ \ \ \ \ \ \ \ \  \ \ \ \ \ \hfill \Box }

\section{Appendix}

\textbf{Proof of \eqref{gammaasym}.}
Recall,
\bqny{
	Z(t,s) =\frac{B(t)+(1-\gamma)B^*(s)}{ct+(1-\gamma)cs+1}, \quad s,t\ge 0,
}	
where $B$ and $B^*$ are independent Brownian motions.
For some positive $\ve$ and large $u$ denote
\begin{equation}\label{ab}
\begin{array}{rcl}
A(\ve) &=&  \Big([0,\IF)\times[0,\IF)\Big) \backslash
\Big([-\ve+t_u,\ve+t_u]\times[0,\ve]\Big),
\\ R(\ve,u) &=&   \Big([-\ve+t_u,\ve+t_u]\times[0,\ve]\Big)\backslash
\Big([-\frac{\ln u}{\sqrt u}+t_u,\frac{\ln u}{\sqrt u}+t_u]
\times[0,\frac{\ln u }{\sqrt u}]\Big).
\end{array}
\end{equation}

We have
\bqn{\label{estim}
\pk{\exists (t,s) \in \mathbb{D}_\delta(u): Z(t,s)>\sqrt u}
&\le&
\pk{\exists (t,s)\in G(\delta/u): Z(t,s)>\sqrt u}
\notag\\&\le&
\pk{\exists (t,s) \in \mathbb{D}_\delta(u): Z(t,s)>\sqrt u}
+
\pk{\exists (t,s) \in A(\ve): Z(t,s)>\sqrt u}
\\&+&
\pk{\exists (t,s): R(\ve,u): Z(t,s)>\sqrt u}
\notag.}
We show that $Z(t,s)$ is a.s bounded for $t,s\ge 0$. According to
Chapter 4, p. 31 in \cite{20lectures} it is equivalent that $Z(t,s)$ is
bounded with positive probability.
We have
\bqny{
\pk{\sup\limits_{t,s\ge 0}Z(t,s)\le 1} &=&
\pk{\text{for all } t,s \ge 0 \ B(t)-ct+(1-\gamma)(B^*(s)-cs) \le 1}
\\ &\ge&
\pk{\text{for all } t,s\ge 0 \ B(t)-ct \le 1/2,(1-\gamma)(B^*(s)-cs) \le 1/2}
\\ &=&
\big(1-\pk{\exists t \ge 0: B(t)-ct>1/2}\big)
\big(1-\pk{\exists s \ge 0: B^*(s)-cs >1/2(1-\gamma)^{-1}}\big)
\\&=&
(1-e^{-c})(1-e^{-c(1-\gamma)^{-1}})>0,
}
where we used \eqref{contclasprob} for the last equation.
Hence by Borell-TIS inequality (see \cite{LifBook})
\bqn{\label{smallborell}
\pk{\exists (t,s) \in A(\ve): Z(t,s)>\sqrt u} = o(\pk{Z(1/c,0)>\sqrt u}), \ \ \ \ u \to \IF
.}

Next we shall prove that
\bqn{\label{smallPiterbarg}
\pk{\exists (t,s) \in R(\ve,u): Z(t,s)>\sqrt u} = o(\pk{Z(1/c,0)>\sqrt u}), \ \ \ \ u \to \IF
.}
If we show that for any $(t,s)\in R(\ve,u)$ and for some positive constant $C$ holds, that
$$\sigma_Z^2(1/c,0)-\sigma_Z^2(t,s) \ge C \frac{\ln^2 u}{u}$$
we can immediately claim \eqref{smallPiterbarg} by Piterbarg's inequality (Proposition 9.2.5 in \cite{20lectures}). Notice that if\\
$i) \ s \notin [0,\frac{\ln u}{\sqrt u}]$, then
\bqny{
\sigma_Z^2(1/c,0)-\sigma_Z^2(t,s) &=&
(A(1/c,0)-A(t,s))+A^*(t,s)
 \\&\ge&
 A^*(t,s) = \frac{s\beta(1-\beta)}{(ct+\beta cs+1)^2}
\ge C\frac{\ln u}{\sqrt u} \ge C \frac{\ln^2 u}{u},
}
hence the claim follows. \\
$ii)$ assume that $s \in [0,\frac{\ln u}{\sqrt u}]$.
Setting
$$L(x) = \frac{x}{(cx+1)^2},$$
we have that  $L(x)$ attains  its unique maxima at point $x=1/c$,
$L'(1/c)=0$ and $L''(1/c)<0$. We have
\bqny{
\sigma_Z^2(1/c,0)-\sigma_Z^2(t,s) =
 \sigma_Z^2(1/c,0)-A(t,s)+A^*(t,s)
\ge
\sigma_Z^2(1/c,0)-A(t,s) = L(1/c)-L(t+\beta s).
}
For all $(t,s)$ such that $(t,s) \in R(\ve,u),
\ s\in [0,\frac{\ln u}{\sqrt u}]$ we have that
$|1/c-(t+\beta s)|\ge C\frac{\ln u}{\sqrt u}$. Hence
$$L(1/c)-L(t+\beta s) \geq C|L''(1/c)|(1/c-(t+\beta s))^2 \geq C \frac{\ln^2 u}{u}$$
and \eqref{smallPiterbarg} holds.
\\
\\
Notice that for some positive constant $C$
\bqny{
\pk{\exists t,s \in \mathbb{D}_\delta(u): Z(t,s)>\sqrt u} \ge
\pk{Z(t_u,0)>\sqrt u} \ge C\pk{Z(1/c,0)>\sqrt u}.
}
Combining the statement above with \eqref{estim},\eqref{smallborell} and \eqref{smallPiterbarg}
we establish \eqref{gammaasym}.
\QED
\\
\\
\textbf{Proof of \eqref{up1}.}
Notice that
$$\pk{\sup\limits_{t\notin [T_u^-,T_u^+]}Z(t)> u} \leq \pk{\exists t,s\in A(\ve):Z(t,s)>\sqrt u}
+\pk{\exists t,s\in R(\ve,u):Z(t,s)>\sqrt u},
$$
where $A(\ve)$ and $R(\ve,u)$ are defined in \eqref{ab}. Hence the claim
follows by \eqref{smallborell} and \eqref{smallPiterbarg}.
\QED
\\
\\
\textbf{Proof of \eqref{intgammasim}.} We shall prove that
\bqn{\label{gammaintsim}
\int\limits_{\R}W_j(x)e^{-\frac{x^2}{2v^2c_i}+x\frac{1+cc_i}{c_i}}dx
=
\int\limits_{-M}^M W_j(x)e^{2cx}dx+\overline{A}_{M,v},
}
where $\overline{A}_{M,v}\to 0$ as $u \to \IF$ and then $M \to \IF$ uniformly for all
$|i|,|j| \leq N_u$.
We have
\bqn{\label{overII} & \ &
\int\limits_{\R}W_j(x)e^{-\frac{x^2}{2v^2c_i}+x\frac{1+cc_i}{c_i}}dx
- \int\limits_{-M}^M W_j(x)e^{2cx}dx
\notag\\&\le&
\int\limits_{|x|>M}W_j(x)e^{-\frac{x^2}{2uc_i}+x\frac{1+cc_i}{c_i}}dx
+|\int\limits_{-M}^M W_j(x)e^{2cx}(
e^{-\frac{x^2}{2uc_i}-\frac{x}{u}\frac{(\theta_u+iS)c}{c_i}}-1)dx|
\notag\\&=:& \overline{I}_1+|\overline{I}_2|.
}

Let $u\ge M^6$. For any integer $|i|,|j|\leq N_u$, all $x\in [-M,M]$ and all $u$ large for some constant $C$ that does not depend on $u$ we have
$$|\frac{(\theta_u+iS)c x}{c_i})|
\leq CM\sqrt u \ln u, \ \ \ \ |\frac{x^2}{2c_i}|\leq CM^2,$$
hence
\bqny{
|\frac{1}{u}(-\frac{x^2}{2c_i} - \frac{(\theta_u+iS)c x}{c_i})| \leq
\frac{C}{u}(M^2+M\sqrt u \ln u)\leq \frac{C}{M^4}+\frac{CM}{u^{2/5}}
\leq \frac{1}{M}.
}
We have by \eqref{finint} and \eqref{epsilon} that   for all $|j|\leq N_u$
$$\int\limits_\R W_j(x)e^{2cx}dx\le C$$
for some constant $C$ that does not depend on $u$
and hence
\bqn{\label{overI2}
|\overline{I}_2| \leq \frac{1}{M}\int\limits_{-M}^M W_j(x)e^{2cx}dx \leq
\frac{1}{M}\int\limits_\R W_j(x)e^{2cx}dx =\frac{C}{M} \to 0, \ \ M \to \IF
.}
Next we have for large $u$
\bqny{
\overline{I}_1&\leq& \int\limits_{|x|>M}
\pk{\exists t\in [0,S],s\ge 0:B(t)-ct+\beta(B^*(s)-cs)>x}e^{x\frac{1+cc_i}{c_i}}dx
\\ &\leq&
\int\limits_{x>M}
\pk{\exists t\in [0,\frac{S}{x}],s\ge 0: Z(t,s)>\sqrt x}e^{x\frac{1+cc_i}{c_i}}dx
+\int\limits_{x<-M}e^{cx}dx
.}
We analyze behavior of $\sigma^2_Z(t,s)$ on the set $\{ (t,s) \in [0,\frac{S}{x}] \times [0,\IF)\}$.
Since
\bqny{
\sigma_Z^2(t,s) = \frac{t+\beta^2s}{(ct+c\beta s+1)^2} \le
\frac{S}{x}+\frac{\beta^2s}{(c\beta s+1)^2} \leq \frac{S}{x}+\frac{\beta}{4c}
}
taking large enough $x$ we can write for any fixed $\ve>0$ that
\bqny{
\sigma_Z^2(t,s) \le \frac{\beta(1+\ve)}{4c}
.}
Hence by Borell-TIS inequality for large $x$
\bqny{
\pk{\exists t\in [0,\frac{S}{x}],s\ge 0: Z(t,s)>\sqrt x}
\leq e^{-x\frac{2c}{\beta (1+2\ve)}}.
}
Choosing $\ve$ such that $\beta (1+2\ve)<1$, uniformly for all $|i|,|j| \le N_u$ we have
with $a = 2c-\frac{2c}{\beta (1+2\ve)}<0$
\bqn{\label{overI1}
\overline{I}_1\leq o(1)+
\int\limits_{|x|>M}
e^{x\frac{1+cc_i}{c_i}-x\frac{2c}{\beta (1+2\ve)}}dx
&=&
o(1)+\int\limits_{|x|>M}
e^{x(2c+o(1)-\frac{2c}{\beta (1+2\ve)})}dx
\notag\\&\le&
o(1) +2\int\limits_{|x|>M}e^{ax}\to 0,  \ \ \ \  \ M \to \IF.
}
Combination of \eqref{overI2} and \eqref{overI1} establishes
\eqref{gammaintsim}. By the monotone convergence theorem \eqref{gammaintsim}
implies \eqref{intgammasim}.
\QED
\\
\\
\textbf{Proof of \eqref{intlim}.}
 We have
\bqny{
\int_{\R}w(x)\omega(j,S,x) dx - \int\limits_{-M}^M w(x) e^{2cx} dx
\le
\int\limits_{|x|>M}w(x) \omega(j,S,x) dx
+  |\int\limits_{|x|<M}w(x) (\omega(j,S,x)
-e^{2cx})|dx
=:  I_1+ |I_2|.
}
Since $W_j(x)\ge w(x)$ we have that ($\overline{I}_1$ and $\overline{I}_2$ are
defined in \eqref{overII})
$$|I_2|\leq |\overline{I}_2|, \ \ \ \ I_1\leq \overline{I}_1$$
implying
$$I_1+I_2 \leq \overline{I}_1+|\overline{I}_2| \to 0 $$
as $u\to \IF$ and then $M \to \IF$ by \eqref{overI2} and \eqref{overI1}. Thus  \eqref{intlim}
is established.
\QED
\\
\\
%Proofs for \eqref{intlimparis} and \eqref{cumlim} are the same since
%$$w(T,x) \leq W_j(x) \ \text{ and } \ w_k(x)\leq W_j(x).$$
\textbf{Proof of \eqref{ratiosum}.}
For any  $j \ge 1$ we have (set $b_j = 2jc^2S$)
\bqny{
\E{\sup\limits_{s \in G^*_j}
e^{\beta(\sqrt 2 B^*(s)-s)}}
&\leq&
\int\limits_\R \pk{\sup\limits_{s\in [2jc^2S,2(j+1)c^2S]}
\beta(\sqrt 2 B(s)-s-\sqrt 2 B(b_j)+b_j)>x+\beta(b_j-\sqrt 2 B(b_j))
}e^xdx
\\ &=&
\frac{1}{\sqrt{2\pi b_j}}
\int\limits_\R\int\limits_\R \pk{\sup\limits_{s\in [0,2c^2S]}
\beta(\sqrt 2 B(s)-s)>x+\beta(b_j-\sqrt 2 y)
}e^xdxe^{-\frac{y^2}{2b_j}}dy
\\ &=&
\frac{1}{\sqrt{2\pi b_j}}
\int\limits_\R\int\limits_\R \pk{\sup\limits_{s\in [0,2c^2S]}
\beta(\sqrt 2 B(s)-s)>x
}e^{x-\beta(b_j-\sqrt 2 y)}e^{-\frac{y^2}{2b_j}}dxdy
\\ &=&
\frac{e^{-\beta b_j}}{\sqrt{2\pi b_j}}
\int\limits_\R\int\limits_\R \pk{\sup\limits_{s\in [0,2c^2S]}
\beta(\sqrt 2 B(s)-s)>x
}e^xe^{-\frac{y^2}{2b_j}+\sqrt 2 \beta y}dxdy
\\ &=&
\frac{e^{-\beta b_j}}{\sqrt{2\pi b_j}}
\int\limits_\R \pk{\sup\limits_{s\in [0,2c^2S]}
\beta(\sqrt 2 B(s)-s)>x
}e^xdx\int\limits_\R
e^{-\frac{y^2}{2b_j}+\sqrt 2 \beta y}dy.
}
For some constant $C$ that does not depend on $S$
we have by \eqref{contclasprob}
\bqny{
\int\limits_\R \pk{\sup\limits_{s\in [0,2c^2S]}
\beta(\sqrt 2 B(s)-s)>x
}e^xdx
\leq
\int\limits_\R \pk{\sup\limits_{s\in [0,\IF)}
 (B(s)-\frac{s}{\beta\sqrt 2})>\frac{x}{\beta\sqrt 2}
}e^xdx
=
1+\int\limits_0^\IF e^{-x/\beta+x} \leq C.
}
Since
$$\int\limits_\R
e^{-\frac{y^2}{2b_j}+\sqrt 2 \beta y}dy = \sqrt{2\pi b_j}e^{\beta^2 b_j},$$
we have for some fixed small enough $\ve$ and large $S$
\bqn{\label{epsilon}
\E{\sup\limits_{s \in G^*_j}
e^{\beta(\sqrt 2 B(s)-s)}} \leq
C e^{-\beta(1-\beta)b_j} \leq e^{-jS\ve}
.}
Thus as $S \to \IF$
\bqny{
\sum\limits_{1\le j\le N_u}
\E{\sup\limits_{s \in G^*_j}
e^{\beta(\sqrt 2 B(s)-s)}}
\leq e^{-\ve S}(1+o(1))
}
establishing the claim.
\QED
\\
\\
{\bf Acknowledgments}: Author would like to thank the anonymous
referee for the useful remarks.
Partial financial supported by  SNSF Grant 200021-175752/1 and PSG 1250
grant (Unil - St. Petersburg) is kindly acknowledged.

\bibliographystyle{ieeetr}
\bibliography{EEEA}{}
\end{document}